\documentclass[twoside,reqno,a4paper,12pt]{amsart}
\usepackage{amsmath,amsthm}
\usepackage{amsmath,amstext,amsfonts}
\usepackage{hyperref}
\usepackage{graphicx}
\usepackage{graphics}
\usepackage{tikz}
\usepackage{etex}

\theoremstyle{plain}

\newtheorem{theorem}{Theorem}
\newtheorem{lemma}{Lemma}
\newtheorem{proposition}{Proposition}
\newtheorem{corollary}{Corollary}
\newtheorem{definition}[theorem]{Definition}
\newtheorem{remark}{Remark}
\newtheorem{example}{Example}[section]

\theoremstyle{remark}

\newcommand{\R}{\mathbb R}

\renewcommand{\qed}{}

\newcommand{\X}{{\mathbb{X}}}
\newcommand{\C}{{\mathbb{C}}}
\newcommand{\N}{{\mathbb{N}}}

\usepackage{floatrow}

\def \to {\rightarrow}

\newcommand{\trip}[1]{{|\kern -1pt\|#1\|\kern -1pt|}}

\begin{document}

\title[Stability Problems in Nonautonomous Linear Differential Equations.]
{Stability Problems in Nonautonomous Linear Differential Equations in Infinite Dimensions.}

\author[Hildebrando M. Rodrigues]
{Hildebrando M. Rodrigues$^{\dag}$}

\thanks{\dag Instituto de Ci\^{e}ncias Matem\'{a}ticas e de
Computa\c{c}ao, Universidade de S\~{a}o Paulo-Campus de S\~{a}o
Carlos, Caixa Postal 668, 13560-970 S\~{a}o Carlos SP, Brazil, e-mail: hmr@icmc.usp.br}

\author[J. Sol\`a-Morales]
{J. Sol\`a-Morales$^{\ddag}$}
\thanks{ \ddag Departament de
Matem\`atiques, Universitat Polit\`ecnica de Catalunya,
Av. Diagonal 647, 08028 Barcelona, Spain, e-mail:
jc.sola-morales@upc.edu}
 \markboth{Joan de Sol\`a-Morales}{
Hildebrando M. Rodrigues}
\thanks{\dag Partially supported by FAPESP Processo 2018/05218-8 }
\thanks{\ddag Partially supported by MINECO grant MTM2017-84214-C2-1-P. Faculty member of the Barcelona Graduate School of Mathematics (BGSMath) and part of the Catalan research group 2017 SGR 01392.}

\author[G. K. Nakassima]
{G. K. Nakassima$^{\dag}$}

\thanks{\dag Instituto de Ci\^{e}ncias Matem\'{a}ticas e de
Computa\c{c}ao, Universidade de S\~{a}o Paulo-Campus de S\~{a}o
Carlos, Caixa Postal 668, 13560-970 S\~{a}o Carlos SP, Brazil, e-mail: hmr@icmc.usp.br, guilherme.nakassima@usp.br}

\maketitle

\begin{center}

{\bf Dedicated to Tom\'as Caraballo for his 60th birthday.}

\end{center}

\begin{abstract} One goal of this paper is to study robustness of stability of nonautonomous linear ordinary differential equations under integrally small perturbations  in an infinite dimensional Banach space.  Some applications are obtained to the case of rapid oscillatory perturbations, with arbitrary small periods, showing that even in this case the stability is robust. 
These results extend to infinite dimensions some results given in Coppel \cite{Coppel2}. Based in Rodrigues \cite{Rodrigues} and in Kloeden \& Rodrigues \cite{KR} we introduce a class of functions that we call Generalized Almost Periodic Functions that extend the usual almost periodic functions and are suitable to deal with oscillatory perturbations.
We also present an infinite dimensional example of the previous results. 
We show in another example that it is possible to stabilize an unstable system using a perturbation with large period and small mean value. Finally, we give an example where we stabilize an unstable linear ODE with small perturbation in infinite dimensions using some ideas developed in Rodrigues \& Sol\`a-Morales \cite{8Ro-Sola}  and in an example of Kakutani, see \cite{Rickart}.

\end{abstract}

\noindent 

\noindent 

\vspace{0.5in}

\section{\bf Introduction} In some papers of some of us we tried to extend or to analyse in infinite dimensions known results for finite dimensional problems. This was the case of Kloeden \& Rodrigues \cite{KR}, Rodrigues \cite{Rodrigues}, Rodrigues \& Ruas \cite{Rodrigues-Ruas}, Rodrigues \& Sol\`a-Morales \cite{Ro-Sola, Ro-Sola1, Ro-Sola2, Ro-Sola3}, Rodrigues, Caraballo \& Gameiro \cite{Rodrigues-Caraballo-Gameiro} and Rodrigues, Teixeira \& Gameiro \cite{Rodrigues-Teixeira-Gameiro}.

Following this philosophy, in this paper we consider the following linear system of  ordinary differential equations in an infinite dimensional Banach space $\X$,
$\dot x=A(t)x$ and a perturbed system $\dot y=A(t)y+B(t)y$, where $A(t)$ and $B(t)$ are continuous in $\R$. We suppose first that for each $t\in \R$, $A(t)$ and $B(t)$ are bounded operators, the first system is asymptotically stable and that $B(t)$ is  integrally small in an arbitrary interval of length bounded by $h>0$. We establish conditions on the smallness of $B(t)$ in such a way that the perturbed system will also be asymptotically stable. This is stablished inTheorem \ref{stabilitybounded}. Then we extend to the case such that $A: {\mathcal D}\to \X$ is ubbounded and generates an $C^0$-semigoup $T(t),\  t\geq 0$. This is stablished inTheorem \ref{stabilityunbounded}.

In Daleckii \& Krein \cite{Daleckii-Krein} page 178 and in Carvalho et al. \cite{Carvalho} similar results are presented about robustness of stability but with the stronger assumption
 $\frac{1}{\tau_0}\int_t^{t+\tau_0}\|B(\tau)\| d\tau\ \leq \delta$, for some $\tau_0>0$,  for every $t\in \R$ for sufficiently small $\delta$. One observes that the smallness condition is imposed with the norm inside the integral and in our case the norm appears outside the integral and this makes a significant difference, as it is shown in Theorem (\ref{stabilitybounded}).

Then we introduce a class of functions that we call Generalised Almost Periodic Functions that contains the usual almost periodic functions. In fact part, of it was introduced in Kloeden \& Rodrigues \cite{KR}, where the authors studied perturbations of an hyperbolic equilibrium.  In the present paper we 
use also the concept of mean value to define the class of Generalized Almost Periodic Functions (GAP).

This new class of functions has some important advantages compared with the almost periodic functions, namely, if we perturb an almost periodic function with a local perturbation in time it will not be almost periodic. Therefore it is not robust with respect to this kind of perturbations. It is also not also robust with respect to some more general perturbations, like chaotic functions. 

As a consequence we study a system of the form
$\dot y=A(t)y+B(\omega t)y$ and prove that if $\omega>0$ is sufficiently large the the stability is preserved. When $B(t)$ is periodic the result says that for sufficiently small periods and large oscillations the stability is preserved.
The function $B(t)$ does not need to be small and so if we have a linear perturbation with large oscillations the stability is preserved. This is shown in Theorem \ref{1stability}.  In the periodic case the perturbation will have very small period. We present an example in the infinite dimensions case, in the space $\ell_2$ where we show that the stability is preserved. These results extend to infinite dimensions some results of Coppel \cite{Coppel2}. 

Then in Theorem (\ref{stabilityunbounded}) we extend the above results  to the case where we have an unbounded infinitesimal generator. 
Henry \cite{Henry} proves similar results with different  applications, but using a different method where he passes from the continuous case to a discrete case and then recover the results for the continuous problem. Our method follows more the method of Coppel \cite{Coppel2} (finite dimension).

In Section \ref{large period} we present a two dimensional example where we show that it is possible to stabilise an unstable system with a periodic perturbation  with large period and small mean value.

Finally in Section \ref{unstable linear ODE} using some ideas developed in Rodrigues \&Sol\`a-Morales \cite{8Ro-Sola} and in an example of Kakutani \cite{Rickart},
we give an example in infinite dimensions where we estabilize  an unstable linear system using a small linear perturbation.

These two last examples seem to be new in the literature, to our knowledge.

\vspace{0.25in}

\maketitle

\section{\bf Robustness of Stability.}\hfill

The next theorem extends to infinite dimensional Banach spaces a result of W. A. Coppel \cite{Coppel2}, Proposition 6, p.6.

\begin{theorem}\label{stabilitybounded}

Let $\X$ be a Banach space and $A,\ B:\ \R\to \ L(\X)$ be continuous functions such that $\|A(t)\|\leq M$ and $\|B(t)\|\leq M$ for every $t\in \R
.$

Consider the equations:
\begin{equation}\label{H}
\dot x = A(t)x
\end{equation}
\begin{equation}\label{P11}
\dot y = A(t)y+B(t)y
\end{equation}
Let $T(t,s)=X(t)X^{-1}(s)$ the evolution operator of (\ref{H}). Suppose that $\|T(t,s)\|\leq Ke^{\alpha (t-s)}$ for $t\geq s$, $t,s\in\ \R$, where $\alpha\in \R$ and $K\geq 1$. 

Let $\delta,\ h$ be two positive numbers.

If $\|\int_{t_1}^{t_2}B(t)dt\|\leq \delta$ for $|t_2-t_1|\leq h$, and $t_1, t_2\in \R$, then the evolution operator $S(t,s)=Y(t)Y^{-1}(s)$ of (\ref{P11}) satisfies the inequality:
$$\|S(t,s)\|\leq (1+\delta)K e^{\beta(t-s)}\  \mbox{for}\  t\geq s, \ t,s\in \R,\ 
\mbox{where}\  \beta=\alpha+3MK\delta+\frac{\log((1+\delta)K)}{h}.$$

If $\alpha$ is negative, h is sufficiently large and $\delta$ sufficiently small in such a way that $\beta <0$ then it follows that system (\ref{P11})
 is asymptotically stable.

\end{theorem}

\noindent{\bf Proof:} By the variation of constants formula $$S(t,s)=T(t,s)+\int _s ^t T(t,u)B(u)S(u,s)du, \ t\geq s.$$

If we let $C(u)=\int_t^uB(\tau)d\tau$
$$\int _s ^t T(t,u)B(u)S(u,s)du=\int _s^tT(t,u)\frac{d}{du}\int_t^uB(\tau)d\tau S(u,s)du=\int _s^tT(t,u)\frac{d}{du}C(u) S(u,s)du$$

Taking derivatives,
$$\frac{d}{du}[T(t,u)C(u)S(u,s)]=$$
$$-T(t,u)A(u)C(u)S(u,s)+T(t,u)B(u)S(u,s)+T(t,u)C(u)(A(u)+B(u))S(u,s)$$

Integrating the above equation, we obtain

$$\int _s ^t\frac{d}{du}[T(t,u)C(u)S(u,s)]du=-\int _s ^tT(t,u)A(u)C(u)S(u,s)du+$$
$$\int _s ^tT(t,u)B(u)S(u,s)du+\int _s ^tT(t,u)C(u)A(u)S(u,s)du+\int _s ^tT(t,u)C(u)B(u)S(u,s)du$$

And so,

$$-T(t,s)C(s)=-\int _s ^tT(t,u)A(u)C(u)S(u,s)du+$$
$$\int _s ^tT(t,u)B(u)S(u,s)du+\int _s ^tT(t,u)C(u)A(u)S(u,s)du+\int _s ^tT(t,u)C(u)B(u)S(u,s)du$$

Therefore,

$$\int _s ^tT(t,u)B(u)S(u,s)du=-T(t,s)C(s)+\int _s ^tT(t,u)A(u)C(u)S(u,s)du$$
$$-\int _s ^tT(t,u)C(u)A(u)S(u,s)du-\int _s ^tT(t,u)C(u)B(u)S(u,s)du.$$

Therefore,
$$S(t,s)=T(t,s)+\int _s ^t T(t,u)B(u)S(u,s)du=$$
$$T(t,s)(I-C(s))+\int _s ^tT(t,u)A(u)C(u)S(u,s)du$$
$$-\int _s ^tT(t,u)C(u)A(u)S(u,s)du-\int _s ^tT(t,u)C(u)B(u)S(u,s)du.$$

We first suppose that $s\leq t\leq s+h$ and estimate $|S(t,s)|$. Let $s\leq u\leq s+h$. Suppose
$$|C(u)|\leq |\int _t^u B(\tau) d\tau|\leq \delta.$$

Therefore,

$$|S(t,s)|\leq K(1+\delta) e^{\alpha(t-s)}+3MK\delta\int_s^t\ e^{-\alpha(t-u)}|S(u,s)|du.$$
and so, using gronwall's inequlity  it follows that  in an arbitrary interval of length $h$, say for $s\leq t\leq s+h$ we have

$$|S(t,s)|\leq K(1+\delta) e^{\alpha (t-s)}e^{3MK\delta (t-s)}=K(1+\delta) e^{(\alpha +3MK\delta) (t-s)}$$

For $t\geq s$ there exists $n\in  \N$,  $n= n(t,s)$ such that $s+nh\leq t\leq s+(n+1)h$ and so 
$$|S(t,s+nh)|\leq K(1+\delta)e^{(\alpha+3MK\delta)(t-s-nh)}.$$

We are going to prove by induction that for $s+nh\leq t\leq s+(n+1)h$
$$|S(t,s)|\leq [K(1+\delta)]^{n+1}e^{(\alpha+3KM\delta)(t-s)}$$

The case $n=0$ has already been proved.

But $S(s+nh,s)=S(s+nh,s+(n-1)h)\cdots S(s+h,s)$ and so
$$|S(s+nh,s)|\leq [K(1+\delta)]^ne^{(\alpha+3KM\delta)nh}$$

Therefore  for $s+nh\leq t\leq s+(n+1)h$
$$|S(t,s)|\leq |S(t,s+nh)||S(s+nh,s)|\leq $$
$${K(1+\delta)}e^{(\alpha +3KM\delta)(t-s-nh)}[K(1+\delta)]^ne^{(\alpha+3KM\delta)nh}= [K(1+\delta)]^{n+1}e^{(\alpha+3KM\delta)(t-s)}$$ 

Therefore for $s+nh\leq t\leq s+(n+1)h$ we have
$$|S(t,s)|\leq[K(1+\delta)]^{n+1}e^{((\alpha+3KM\delta)(t-s)}.$$

Let $\gamma \doteq \frac{\ln((1+\delta) K)}{h}$. Since $t\geq s+nh$, we have
$$[(1+\delta)K]^n=e^{\gamma nh}\leq e^{\gamma(t-s)}.$$

Therefore,
$$|S(t,s)|\leq K(1+\delta)e^{(\alpha +3KM\delta+\frac{\ln((1+\delta) K)}{h})(t-s)}\qed$$

\section{\bf The space of generalised almost periodic functions}\label{funsp}

Let $(\mathbb{X},|\cdot|)$ be a Banach space and recall the definition of an almost periodic function \cite{Fink}.
\begin{definition}\label{ap}
A continuous function $f:\mathbb R \to \mathbb X$ is said to be almost periodic
if for every sequence $(\alpha^{\prime} _n)$ there exists a subsequence $(\alpha _n)$ such that the
$\lim_{n\to \infty}f(t+\alpha _n)$ exists uniformly in $\R$.
\end{definition}

Now let  $BUC(\mathbb{R},L(\mathbb X)$ denote the space of bounded and uniformly continuous  functions  $A$  $:$ $\mathbb{R}$ $\to$ $L(\mathbb X)$, which  is a Banach space with the supremum norm 
$\|A\| \doteq \sup_{t\in \R}|A(t)|$, and define
$$
\mathcal{F} \doteq \left\{A \in BUC(\mathbb{R},L(\mathbb X))\ : \,A \,\,
\mbox{is uniformly continuous with precompact range}  \,
\mathcal{R}(A)  \right\}.
$$
The class $\mathcal{F}$ is quite large and includes both periodic and almost periodic functions as well as
other nonrecurrent functions.
\begin{proposition}
Let $A(t)\in L(\X)$ be almost periodic. Then $A\in {\mathcal F}$.
 \end{proposition}
\noindent {\bf Proof:} The proof is trivial.
\begin{theorem}\label{Fbanach}
$\mathcal{F}$ is a closed subspace of $BUC(\mathbb{R},L(\mathbb X))$ and
hence a Banach space.
\end{theorem}
{\bf Proof:} This proof can be found in Kloeden-Rodrigues \cite{KR}.

\begin{lemma}
Let  $\sup_{t\in \R}|A(t)|\leq M$, If there exists $\lim_{T\to \infty}\frac{1}{T}\int_{a}^{a+T}A(t)dt$ for some  $a\in \R$
then it is independent of $a$.
\end{lemma}

\noindent {\bf Proof:} Let $a\in \R$. $$|\frac{1}{T}\int_{a}^{a+T}A(t)dt-\frac{1}{T}\int_{0}^{T}A(t)dt|=|\frac{1}{T}[\int_{a}^{a+T}A(t)dt-\int_{0}^{T}A(t)dt]|$$
$$|\frac{1}{T}[\int_{a}^{0}A(t)dt+ \int_{T}^{a+T}A(t)dt]|\leq \frac{2M|a|}{T}\to 0,\ \mbox{as}\ T\to \infty.$$ \qed

Then we define:
\begin{definition}
We say that $A\in {\mathcal {F}}$ is a generalized almost periodic function if there exists the limit
$\lim_{T\to \infty}\frac{1}{T}\int_{a}^{a+T}A(t)dt$ in  $L(\mathbb X)$, that is, there exists $\mathbf A \in L(\mathbb X)$ such that, given $\varepsilon >0$ there exists $T_0 =T_0(\varepsilon)>0$ such that $|\frac{1}{T}\int_{a}^{a+T}A(t)dt-\mathbf A|<\varepsilon$ for every $T\geq T_0$ uniformly with respect do $a\in \mathbb R$.

\end{definition} 
\begin{definition}
We define the class of generalized almost periodic functions as
$$\mathcal{GAP}= \{A\in {\mathcal {F}}:\ A\ \mbox{is a generalized almost periodic function}\}$$
\end{definition}
\begin{lemma}
$\mathcal{GAP}$ is a closed subspace of  ${\mathcal {F}}$.
\end{lemma}
{\bf Proof:} Let $A_n \in  \mathcal{GAP}$, $A_n \to A$ in $ {\mathfrak {F}}$.
We must prove that $A\in \mathcal{GAP}$. Given $\varepsilon >0$ there exists $n_0=n_0(\varepsilon)$ such that
$\|A_{n_0}-A\|=\sup_{t\in \mathbb R}|A_{n_0}(t)-A(t)| <\varepsilon$.

Since there exists the $\lim_{T\to \infty}\frac{1}{T}\int_{a}^{a+T}A_{n_0}(t)dt={\mathbf A}_{n_0}$,
there exists $T_0=T_0(\varepsilon)$ such that
$$T_1,\ T_2>T_0 \ \Rightarrow |\frac{1}{T_2}\int_{a}^{a+T_2}A_{n_0}(t)dt-\frac{1}{T_1}\int_{a}^{a+T_1}A_{n_0}(t)dt|<\varepsilon,\ \forall a\in \mathbb R$$

Then
$$T_1,\ T_2>T_0 \ \Rightarrow |\frac{1}{T_2}\int_{a}^{a+T_2}A(t)dt-\frac{1}{T_1}\int_{a}^{a+T_1}A(t)dt|\leq$$
$$|\frac{1}{T_2}\int_{a}^{a+T_2}A(t)dt-\frac{1}{T_2}\int_{a}^{a+T_2}A_{n_0}(t)dt|+|\frac{1}{T_2}\int_{a}^{a+T_2}A_{n_0}(t)dt-\frac{1}{T_1}\int_{a}^{a+T_1}A_{n_0}(t)dt|+$$
$$|\frac{1}{T_1}\int_{a}^{a+T_1}A_{n_0}(t)dt-\frac{1}{T_1}\int_{a}^{a+T_1}A(t)dt|\leq$$
$$\frac{1}{T_2}\int_{a}^{a+T_2}|A(t)-A_{n_0}(t)|dt+|\frac{1}{T_2}\int_{a}^{a+T_2}A_{n_0}(t)dt-\frac{1}{T_1}\int_{a}^{a+T_1}A_{n_0}(t)dt|+$$
$$|\frac{1}{T_1}\int_{a}^{a+T_1}|A(t)-A_{n_0}(t)|dt\leq 3\varepsilon$$

Using Cauchy Criterion we conclude that there exists
$$\lim_{T\to \infty}\frac{1}{T}\int_{a}^{a+T}A(t)dt={\mathbf A} \in L(\mathbb X), \ \forall a\in \mathbb R$$

This implies that $A\in \mathcal{GAP}$.  \qed
\begin{definition}
For  $A\in \mathcal{GAP}$ we define the mean value of $A$ as:
 $${\mathcal M}(A)\doteq \lim_{T\to \infty}\frac{1}{T}\int_{a}^{a+T}A(t)dt \in L(\X).$$

\end{definition}

\begin{lemma} The function ${\mathcal M}:\mathcal{GAP}\to L(\X)$ is an uniformly continuous function.
\end{lemma}
\noindent {\bf Proof:} Let $A,B\in \mathcal{GAP}$. Then
$$|{\mathcal M}(A)-{\mathcal M}(B)|=| \lim_{T\to \infty}\frac{1}{T}\int_{a}^{a+T}A(t)dt- \lim_{T\to \infty}\frac{1}{T}\int_{a}^{a+T}B(t)dt|=$$
$$| \lim_{T\to \infty}\frac{1}{T}\int_{a}^{a+T}[A(t)-B(t)]dt|\leq \sup_{t\in \mathbb R}|A(t)-B(t)|=\|A-B\|.$$ \qed

 Let $\mathcal O=\{A\in\ \mathcal{GAP}:{\mathcal M}(A)=\lim _{T\to \infty}\frac{1}{T}\int_{a}^{a+T}A(t)dt=0, \ \forall \ a\in \mathbb R\}$
\begin{corollary}
$\mathcal O$ is a closed subspace of $\mathcal{GAP}$.
\end{corollary}
\noindent {\bf Proof:} Since ${\mathcal M}(A)$ is a continuous function, the set $\mathcal O={\mathcal M}^{-1}\{\{0\}\}$ is closed set.
\begin{corollary}
Any function $A\in \mathcal{GAP}$ can be written as $A=A_0+B$, where $A_0={\mathcal M}(A)$ and $B\in \mathcal O$.
\end{corollary}

The next theorem shows that stability is preserved if the linear perturbation has sufficiently large frequency:

\begin{theorem}\label{1stability}

Let $A,\ B:\ \R\to \ L(\X)$ be continuous functions such that $\|A(t)\|\leq M$ and $\|B(t)\|\leq M$ for every $t\in \R.$
Suppose that $B(t)$ is a generalized almost periodic function with mean value zero ($\mathcal{GAP}$).
Consider the equations:
\begin{equation}\label{Hw}
\dot x = A(t)x
\end{equation}
\begin{equation}\label{Pw}
\dot x = A(t)x+B(\omega t)x
\end{equation}
Let $T(t,s)$ the evolution operator of (\ref{Hw}). Suppose that $\|T(t,s)\|\leq Ke^{-\alpha (t-s)}$ for $t\geq s$, $t,s\in\ \R$, where $\alpha>0$ and $K>1$.
Then there exists $\tilde K$ and $\omega _0>0$ such that for $\omega>\omega _0$ 
$$|S_{\omega}(t,s)|\leq \tilde K e^{\frac{-\alpha}{2}(t-s)},\ t\geq s,$$
where $S_{\omega}(t,s)$ indicates the evolution operator of (\ref{Pw}).
\end{theorem}

\noindent {\bf Proof:}

We are going to show that for any $h>0$, $\delta >0$ there exists $\omega _0=\omega _0(h,\delta)>0$
such that if $\omega >\omega _0$ then
$$|\int_{t_1}^{t_2}B(\omega t)dt|\leq \delta\ \mbox{for}\ |t_2-t_1|\leq h.$$

Let us consider first the case $|t_2-t_1|\leq \frac{\delta}{M}$. Since $|B(t)|\leq M$ for every $t\in \R$, we have
$$|\int_{t_1}^{t_2}B(\omega t)dt|\leq |\int_{t_1}^{t_2}|B(\omega t)|dt|\leq M|t_2-t_1|\leq M\frac{\delta}{M}=\delta.$$
To complete the proof we consider now the case $h\geq |t_2-t_1|\geq \frac{\delta}{M}$.

Since $B(t)$ has mean value zero, there exists $T_0=T_0(\frac{\delta}{h})>0$ such that 
$$T\geq T_0\ \Rightarrow |\frac{1}{T}\int _s^{s+T}B(t)dt|\leq \frac{\delta}{h}\ \mbox{for all } s\in \R.$$

By a change of variables,

$$\int_{t_1}^{t_2}B(\omega t)dt=\frac{1}{\omega}\int_{\omega t_1}^{\omega t_2}B(u)du$$
and so for $\frac{\delta}{M}\leq |t_2-t_1| \leq h$
$$|\int_{t_1}^{t_2}B(\omega t)dt|=\frac{1}{\omega |t_2-t_1|}|\int_{\omega t_1}^{\omega t_2}B(u)du||t_2-t_1|\leq \frac{1}{ |\omega t_2-\omega t_1|}|\int_{\omega t_1}^{\omega t_2}B(u)du|\ h.$$

If we take $\omega _0\doteq \frac{MT_0}{\delta}$ we have for $\omega \geq \omega _0$
$$|\omega t_2-\omega t_1|\geq \omega _0|t_2-t_1|\geq\frac{MT_0}{\delta}\frac{\delta}{M}=T_0$$

Therefore,
$$|\int_{t_1}^{t_2}B(\omega t)dt|=\frac{1}{ |\omega t_2-\omega t_1|}|\int_{\omega t_1}^{\omega t_2}B(u)du|\ h\leq \frac{\delta}{h}\ h=\delta$$  

The result follows from Theorem \ref{1stability} for $\delta$ sufficiently small. \qed

Consider now $A\in \mathcal {GAP}$. Then we have $A(t)=A_0+B(t)$, where $A_0={\mathcal M}(A)$ and 
${\mathcal M}(B)=0$. We suppose that $|A_0|\leq M$ and $|B(t)|\leq M$ for every $t\in \R$. Consider the equations:
\begin{equation}\label{H1w}
\dot x = A_0x
\end{equation}
\begin{equation}\label{P1w}
\dot x = A_0x+B(\omega t)x
\end{equation}

Let $T(t)\doteq e^{A_0 t}$ be the semigroup generated by (\ref{H1w}) and $S_{\omega}(t,s)$ be the evolution operator of (\ref{P1w}).

As a consequence of Theorem \ref{1stability} it follows that if $\sigma (A_0)\subset \{\lambda \in \C: Re(\lambda)<-\alpha\}$ we willl have:
\begin{corollary}
Let  $|T(t)|\leq Ke^{-\alpha t}$ for $t\geq 0$, $K\geq 1$. Then there exists $\widetilde\alpha<\alpha$,
$\widetilde K>K$, $\omega_0 >0$, such that for $\omega>\omega_0$ we have 
$$S_{\omega} (t,s)\leq \tilde{K}e^{-\widetilde{\alpha}(t-s)}, \forall t\geq s.$$

\end{corollary}

\section{\bf An infinite dimensional example}\label{inf dim example}

In this section we will construct a true infinite dimension example to apply the results of the previous section. We are going to use
some results of the paper Rodrigues and Sol\`a-Morales \cite{Ro-Sola2}. Consider the space $\X=\ell _2$. We consider the operator
$J\in {\mathcal L}(X)$ given by the infinite dimensional Jordan matrix:
\begin{equation}\label{Jordan}
J:=\begin{pmatrix}
0&0&0&\cdots \\
1&0&0&\cdots \\
0&1&0&\cdots \\
0&0&1&\cdots \\
.&.&.&\cdots \\
.&.&.&\cdots 
\end{pmatrix}
\end{equation}

As it is proved in Rodrigues and Sol\`a-Morales \cite{Ro-Sola2} the spectrum of $J$ is the closed unity circle of the complex plane.
Now we take $0<a<1$ and we define the operator:
\begin{equation}\label{Loperator}
L:=\begin{pmatrix}
a&0\\
0&\nu J+aI
\end{pmatrix}=aI+\nu \begin{pmatrix} 0&0\\
0&J \end{pmatrix}=a\left (I-\left (\frac{-\nu}{a}\right )\begin{pmatrix} 0&0\\
0&J \end{pmatrix}\right )
\end{equation}

If we let 
$$D=\left (\frac{-\nu}{a}\right )\begin{pmatrix} 0&0\\
0&J \end{pmatrix}$$
we have that
$$L=a(I-D).$$

 From the same paper above it follows that the spectrum of $L$ is the closed disc $B_{\nu}(a)$ with center in $a$ and radius $\nu$. Then we take
 $0<\nu<\min\{a,1-a \}$

Then we let $A:=\log L= (\log a)I+\log(I-D)$.

\begin{figure}[!h]
\begin{floatrow}
\includegraphics[height=0.3\textheight]{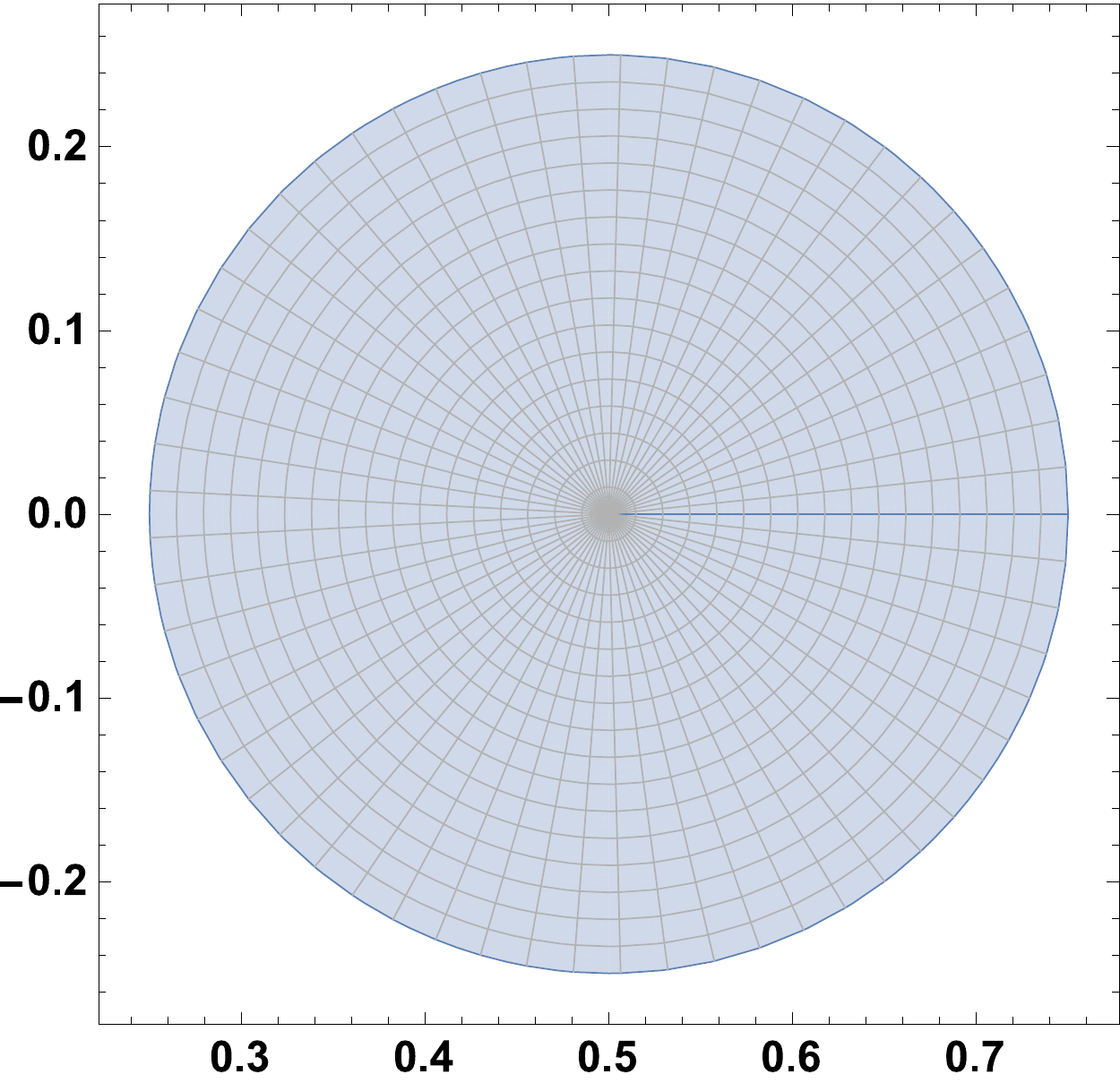}
\hfill
\includegraphics[height=0.3\textheight]{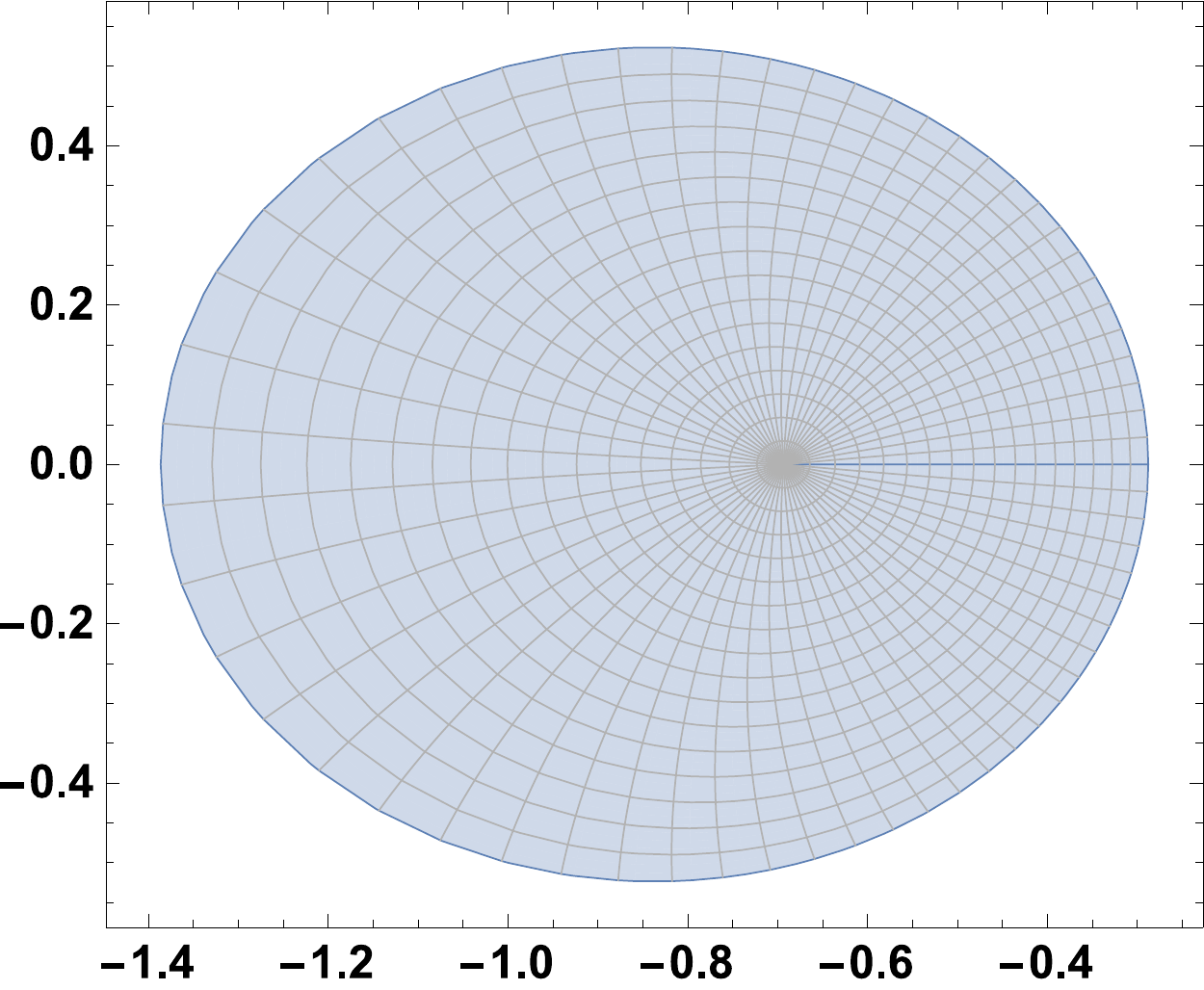}
\caption{{\bf Left}: The spectrum of $L$ given by $\sigma (L) = B_{\nu}(a)$.
{\bf Right}: The spectrum of $A$ given by $\sigma (A) = \log (\sigma (L))$,
with $a = 1/2$ and  $\nu=1/4$.}
\label{fig:spectrum_A}
\end{floatrow}
\end{figure}


But
$$\log(I-D)=-(D+\frac{D^2}{2}+\cdots \frac{D^n}{n}\cdots).$$
Therefore 
$$\|\log(I-D)\|\leq \frac{\nu}{a}+\frac{( \frac{\nu}{a})^2}{2}+\cdots \frac{( \frac{\nu}{a})^n}{n}+\cdots =-\log (1-\frac{\nu}{a})$$

Let $\nu>0$ sufficiently small such that  $0<-\log (1-\frac{\nu}{a})<\frac{a}{2}$.

Then it follows that
$$\|e^{At}\|\leq e^{(-at-\log (1-\frac{\nu}{a})t)}\leq e^{-\frac{a}{2}t},\ \ t\geq 0$$

In the space $\X=\ell _2$. We consider the operator
$A\in {\mathcal L}(X)$ given above.

\begin{corollary}\label{B  GAP}
Consider now the systems:
\begin{equation}\label{H2w}
\dot x= Ax
\end{equation}
\begin{equation}\label{P2w}
\dot y= Ay+B(\omega t)y
\end{equation}
where $B\in {\mathcal GAP}$ with meanvalue zero. Let $M>0$ be such that $|A|\leq M$ and $\sup_{t\in \R}|B(t)|\leq M$.

Let $S_{\omega}(t,s)=Y(t)Y^{-1}(s)$ be the evolution operator associated to to system (\ref{P2w}), where $Y(t)$ is the solution with initial condition $Y(0)=I$, where $I$ indicates the Identity operator.

Then there exists $\tilde K$, $\tilde \alpha$ and $\omega _0>0$ such that for $\omega>\omega _0$ 
$$|S_{\omega}(t,s)|\leq \tilde K e^{- {\tilde\alpha}(t-s)},\ t\geq s,$$

\end{corollary}

{\bf Proof:} Follows from  Theorem \ref{1stability} . \qed

Next we will present a simple example where the perturbation $B(t)$ belongs to  ${\mathcal GAP} $ but it is not almost periodic.

\begin{example}\label{infinite example}
Let $b:\R\to \R$ be uniformly continuous, bounded with mean value zero. Let 
\begin{equation}
B(t)\doteq \begin{pmatrix}
0&0&0&\cdots \\
b(t)&0&0&\cdots \\
0&b(t)&0&\cdots \\
0&0&b(t)&\cdots \\
.&.&.&\cdots \\
.&.&.&\cdots 
\end{pmatrix}=b(t)
\begin{pmatrix}
0&0&0&\cdots \\
1&0&0&\cdots \\
0&1&0&\cdots \\
0&0&1&\cdots \\
.&.&.&\cdots \\
.&.&.&\cdots 
\end{pmatrix}
\end{equation}

Then $B\in \mathcal{GAP}$ and has mean value zero.
Let $d(t)\doteq \sqrt{1-t^2}$ if $-1\leq t\leq 1$, $d(t)=0$ if $t\in (\-\infty,0)\cup (1,\infty)$. In the special case that we take $b(t)\doteq d(t)+\cos t$, $B(t)$ is not almost periodic.

Therefore we can apply Corollary \ref{B GAP} if we take $b(\omega t)=d(\omega t)+\cos(\omega t)$ and then we can take $B(\omega t)$.
\end{example}

\section{\bf A case where the infinitesimal generator is unbounded.}

Consider the equations:
\begin{equation}\label{HA}
\dot x = Ax
\end{equation}
\begin{equation}\label{PB}
\dot y = Ay+B(t)y
\end{equation}
We suppose that $\mathcal D$ is dense in $\X$ and $A:{\mathcal D}\to \X$ is the infinitesimal generator of a ${\mathcal C}_0$ semigroup $T(t)$, such that $|T(t)|\leq Ke^{\alpha t},\ t\geq 0,\  K\geq 1, \alpha \in \R$.

Now we will analyse some  smallness conditions on the perturbation $B(t)$, such that the equation \ref{PB} is also asymptotically stable in the case $\alpha<0$. The case when $B(t)$ is uniformly small is studied in Kloeden-Rodrigues \cite{KR} without leaving the continuous case. Similar results are obtained by Carvalho et all \cite{Carvalho}, but they first find the result for the discrete case.

Similar results to the next theorem are treated by Carvalho et all \cite{Carvalho} and Dalekii-Krein \cite{Daleckii-Krein} but they use the stronger assumption that $\int _{\tau}^t|B(t)|$ is small, with the norm inside the integral and in the first one they prove via a discretiztion method. Similar results are obtained by Henry \cite{Henry} in Thorem7.6.11, pag. 238, where he also consider first the discrete case, and requires that $B(t)$ is uniformly small and integrally small.

Our result is an extension of a classical result of Coppel \cite{Coppel2} for the infinite dimensional case, and $A$ being an unbounded operator.

We will follow the steps of Theorem \ref{stabilitybounded} where we imposed that $|B(t)|\leq M$ for every $t\in \R$ and that $|\int_t^uB(\tau)d\tau|\leq \delta$ for $t\leq u \leq t+h$. We also assume that the range of $B(t)$ is contained in the domain of $A$.  

\begin{theorem}\label{stabilityunbounded}

We assume besides the above assumptions on $A$ and $T(t)$, $B:\ \R\to \ L(\X)$ is a continuous function and such that for each $t\in \R$ $AB(t)$ is a bounded operator and $B(t)A$ can be extended to the whole space as a bounded operator. We suppose also that $AB(t)$ and $B(t)A$ are bounded for every  $t\in \R$.
 For each $t\in \R$ let $C_t(u) \doteq \int_t^uB(\tau)d\tau$, for $|t-u|\leq h$, where $h$ is a positive real number. We suppose that   there are positive numbers $M$ and $\delta$ such that
 $$|C_t(u)B(u)|\leq M\delta,\  |C_t(u)A|\leq M\delta, \mbox{and}\ |AC_t(u)|\leq M\delta,\  \mbox{for}\  |u-t|\leq h.$$
 
Let $S(t,s)$ be the evolution operator associated to system \ref{PB}. Then
$$\|S(t,s)\|\leq (1+\delta)K e^{\beta(t-s)}\  \mbox{for}\  t\geq s, \ t,s\in \R,\ 
\mbox{where}\  \beta=\alpha+3MK\delta+\frac{\log((1+\delta)K)}{h}.$$
If $\alpha$ is negative, h is sufficiently large and $\delta$ sufficiently small in such a way that $\beta <0$ then it follows that system (\ref{PB})
 is asymptotically stable.

\end{theorem}

\noindent{\bf Proof:} The proof follows the ideas of (\ref{stabilitybounded}). By the variation of constants formula $$S(t,s)=T(t-s)+\int _s ^t T(t-u)B(u)S(u,s)du, \ t\geq s.$$
$$\int _s ^t T(t-u)B(u)S(u,s)du=\int _s^tT(t-u)\frac{d}{du}\int_t^uB(\tau)d\tau S(u,s)du=\int _s^tT(t-u)\frac{d}{du}C_t(u) S(u,s)du$$

Taking derivatives,
$$\frac{d}{du}[T(t-u)C_t(u)S(u,s)]=$$
$$-T(t-u)AC_t(u)S(u,s)+T(t-u)B(u)S(u,s)+T(t-u)C_t(u)(A+B(u))S(u,s)$$

Integrating the above equation, we obtain

$$\int _s ^t\frac{d}{du}[T(t-u)C_t(u)S(u,s)]du=-\int _s ^tT(t-u)AC_t(u)S(u,s)du+$$
$$\int _s ^tT(t-u)B(u)S(u,s)du+\int _s ^tT(t-u)C_t(u)AS(u,s)du+\int _s ^tT(t-u)C_t(u)B(u)S(u,s)du$$

And so,

$$-T(t-s)C_t(s)=-\int _s ^tT(t-u)AC_t(u)S(u,s)du+$$
$$\int _s ^tT(t-u)B(u)S(u,s)du+\int _s ^tT(t-u)C_t(u)AS(u,s)du+\int _s ^tT(t-u)C_t(u)B(u)S(u,s)du$$

Therefore,

$$\int _s ^tT(t-u)B(u)S(u,s)du=-T(t-s)C_t(s)+\int _s ^tT(t-u)AC_t(u)S(u,s)du$$
$$-\int _s ^tT(t-u)C(_tu)AS(u,s)du-\int _s ^tT(t-u)C_t(u)B(u)S(u,s)du.$$

Therefore,
$$S(t,s)=T(t-s)+\int _s ^t T(t-u)B(u)S(u,s)du=$$
$$T(t-s)(I-C_t(s))+\int _s ^tT(t-u)AC_t(u)S(u,s)du$$
$$-\int _s ^tT(t-u)C_t(u)AS(u,s)du-\int _s ^tT(t-u)C_t(u)B(u)S(u,s)du.$$

We first suppose that $s\leq t\leq s+h$ and estimate $|S(t,s)|$. 

If $0\leq |u-t|\leq h$ then 

$$|C_t(u)B(u)|\leq |\int _t^u B(\tau) d\tau B(u)|\leq M\delta,\ |C_t(u)A|\leq M\delta\  \mbox{and}\  |AC_t(u)|\leq M\delta $$

Therefore,
$$|S(t,s)|\leq K(1+\delta) e^{\alpha(t-s)}+3MK\delta\int_s^t\ e^{-\alpha(t-u)}|S(u,s)|du.$$
and so using Gronwall's inequality it follows that  in an arbitrary interval of length $h$, say for $s\leq t\leq s+h$ we have

$$|S(t,s)|\leq K(1+\delta) e^{\alpha (t-s)}e^{3MK\delta (t-s)}=K(1+\delta) e^{(\alpha +3MK\delta) (t-s)}$$

For $t\geq s$ there exists $n\in  \N$, $n= n(t,s)$ such that $s+nh\leq t\leq s+(n+1)h$ and so 
$$|S(t,s+nh)|\leq K(1+\delta)e^{(\alpha+3MK\delta)(t-s-nh)}.$$

We are going to prove by induction that for $s+nh\leq t\leq s+(n+1)h$
$$|S(t,s)|\leq [K(1+\delta)]^{n+1}e^{(\alpha+3KM\delta)(t-s)}$$

The case $n=0$ has already been proved.

But $S(s+nh,s)=S(s+nh,s+(n-1)h)\cdots S(s+h,s)$ and so
$$|S(s+nh,s)|\leq [K(1+\delta)]^ne^{(\alpha+3KM\delta)nh}$$

Therefore  for $s+nh\leq t\leq s+(n+1)h$
$$|S(t,s)|\leq |S(t,s+nh)||S(s+nh,s)|\leq $$
$${K(1+\delta)}e^{(-\alpha -3KM\delta)(t-s-nh)}[K(1+\delta)]^ne^{(\alpha+3KM\delta)nh}= [K(1+\delta)]^{n+1}e^{(\alpha+3KM\delta)(t-s)}$$ 

Therefore for $s+nh\leq t\leq s+(n+1)h$ we have
$$|S(t,s)|\leq[K(1+\delta)]^{n+1}e^{(\alpha+3KM\delta)(t-s)}.$$

Let $\gamma \doteq \frac{\ln((1+\delta) K)}{h}$. Since $t\geq s+nh$, we have
$$[(1+\delta)K]^n=e^{\gamma nh}\leq e^{\gamma(t-s)}.$$

Therefore,
$$|S(t,s)|\leq K(1+\delta)e^{(\alpha -+3KM\delta+\frac{\ln((1+\delta) K)}{h})(t-s)}\qed$$

\section{\bf Applications}
Consider the following result from Henry  \cite{Henry} pg. 30.
\begin{theorem}
Suppose A is a closed operator in the Banach space $\X$ and suppose that $\sigma _1$ is a bounded spectral set of $A$,\ and $\sigma _2=\sigma(A)- \sigma _1$ so $\sigma _2 \cup\{\infty \}$ is another spectral set. Let $E_1$, $E_2$ be the projections associated with these spectral sets, and $X_ j=E_j (\X)$, $j=1,2$. Then $\X=X_1 \oplus X_2$, the $X_j$ are invariant under A, and if $A_j$ is the restriction of $A$ to $X_j$, then 
$$A_ 1: X_1\to X_1 \ \mbox{is bounded}, \sigma(A_1)=\sigma _1,\ \ {\mathcal D}(A_2)={\mathcal D}(A)\cap X_2 \ \mbox{and}\ \sigma(A_2)=\sigma _2.$$
\end{theorem}

\begin{theorem} \label{stabilityunbounded1}
Let $h$ and $\delta$ be positive real numbers.

Suppose that  $A:{\mathcal D}(A)\subset \X\to \X$ a generator of a $C_0$-semigroup $T(t),\ t\geq 0$, $B(t)\in L(\X)$ and
$|B(t)|\leq M$  for every $t\in \R$. Suppose we can decompose $\sigma(A)\doteq \sigma _1\cup \sigma _2$, where $\sigma _1$ is a bounded spectral set and $\sigma _2=\sigma(A)- \sigma _1$ so $\sigma _2 \cup\{\infty \}$ is another spectral set. Suppose there is a smooth curve $\Gamma$, oriented positively, that contains $\sigma_1$ in its interior and $\sigma_2$ is in the exterior of $\Gamma$. 
Consider the projection $P_1\doteq \frac{-1}{2 \pi i}\oint _{\Gamma} (\lambda -A)^{-1}d\lambda$ that projects $\X$ in the subspace $X_1$ associated to the spectral set $\sigma _1$. Let $P_2\doteq I-P_1$. $|T(t)P_1|\leq Ke^{-\alpha t}$ and $|T(t)P_2|\leq K e^{-\mu t}$, for $t\geq 0$, where $\mu>\alpha$. Then $AP_1$ is a bounded operator and $P_1A=AP_1$ and so $P_1A$ is also a bounded operator.

The above decomposition is chosen in such a way that $|P_2B(t)|\leq M\delta$ for every $t\in \R$ .

In analogy with the bounded case if $C_t(u) \doteq \int_t^uB(\tau)d\tau$,  we suppose that  
\begin{equation}
|P_1C_t(u)B|\leq M\delta,\ |P_1AC_t(u)|\leq M\delta\  \mbox{and}\  |P_1C_t(u)A|\leq M\delta,\  \mbox{for}\  t\leq u\leq t+h.\end{equation}\label{smallness conditions}
  Consider the equations:
\begin{equation}\label{H2}
\dot x=Ax
\end{equation}
\begin{equation}\label{P2}
\dot y=Ay+B(t)y
\end{equation}

If the above assumptions are satisfied if $\delta$ is sufficiently small, $h$ is sufficiently large and
(\ref{H2}) is asymptotically stable then system (\ref{P2}) is also asymptotically stable.
\end{theorem}
{\bf Proof:}  The proof follows the ideas of Theorem \ref{stabilityunbounded}.

\begin{remark}
The decomposition $\sigma(A)=\sigma_1\cup \sigma _2$ and the smallness conditions (\ref{smallness conditions}) are satisfied if 
$A$ is at least a sectorial operator  and if $B(t)$ comutes with $P_1$.
\end{remark}

\section{\bf Stabilising unstable systems under small periodic perturbation,\\
 with large period.} \label{large period}\hfill
The next example is in $\R^2$ and it shows that it is possible stabilise an unstable system under a small (in mean value)
 periodic perturbation.
 
 Let $0 < \alpha < \beta$ and $\delta < T$. Let 
 $$
    A \doteq \begin{pmatrix}
            \alpha  & 0\\
            0       & -\beta
        \end{pmatrix},\ \ \ 
    R  \doteq \begin{pmatrix}
            0                    & \frac{\pi}{2\delta}\\[2pt]
            -\frac{\pi}{2\delta} & 0
        \end{pmatrix}
$$
 
 Let $D(t)$ the $T$-periodic operator given by
 \begin{equation}\label{D perturb}
    D(t) =-A +R,\ T-\delta \leq t <T,\ D(t)=0\  t\in \R-[T-\delta, T).
   \end{equation}
 
 Consider the systems:
 \begin{equation}\label{Hper}
 \dot x=Ax
 \end{equation}
 \begin{equation}\label{Pper}
 \dot y=Ay+D(t)y
 \end{equation}

 First we observe that $\lim_{T\to \infty} \frac{1}{T}\int_0^T D(s)ds=0$, that is $B(t)$ has zero mean value, but has large period.
 
 Nest we are going to prove, using Floquet Theorem that system (\ref{Pper}) is uniformly asymptotically stable.
 
 For the sistem $\dot x =A(t)x$, where $A(t)$ is continuous and $T$-periodic, will use Floquet's Theorem even if , $A(t)$ is not continuous, according to the comment in \cite{Hale} page 118.
 
 Consider the matrix solution $X(t)$ of (\ref{Hper}) such that  $X(0)=I$ the identity matrix. Then it is given by
 
  \begin{equation*}
        X(t) = e^{At} = \begin{pmatrix}
                            e^{\alpha t}& 0\\
                            0           & e^{-\beta t}
                        \end{pmatrix}
    \end{equation*}
   
 If we let $R\doteq \begin{pmatrix}
                            0& \frac{\pi}{2\delta}\\
                            -\frac{\pi}{2\delta}        & 0
                        \end{pmatrix}
$  then we have  the rotation matrix:
   
    \begin{equation*}
        e^{Rt} = \begin{pmatrix}
                            \cos(\frac{\pi t}{2\delta})& \sin(\frac{\pi t}{2\delta})\\
                            -\sin(\frac{\pi t}{2\delta})          &  \cos(\frac{\pi t}{2\delta})
                        \end{pmatrix}
    \end{equation*}
Since $X(T-\delta)=e^{A(t-\delta)}=\begin{pmatrix}
        e^{\alpha(T-\delta)}& 0 \\
        0                   & e^{-\beta(T-\delta)}
    \end{pmatrix},$
the fundamental matrix $Y(t)$ of $\dot y=(A+D(t))y$, such that $Y(0)=I$ will be given by 
$$Y(t)=e^{At}\ \mbox{for} \ 0\leq t<T-\delta,\ \ Y(t)=e^{R(t-(T-\delta))}e^{A(T-\delta)}=e^{R(t-T)}e^{R\delta}e^{A(T-\delta)},$$
$\mbox{for} \ T-\delta\leq t<T.$

Then the monodromy matrix will be
$$Y(T)=e^{R\delta}e^{A(T-\delta)}= \begin{pmatrix}
                            \cos(\frac{\pi \delta}{2\delta})& \sin(\frac{\pi \delta}{2\delta})\\
                            -\sin(\frac{\pi \delta}{2\delta}) &  \cos(\frac{\pi \delta}{2\delta})
                        \end{pmatrix}\begin{pmatrix}
        e^{\alpha(T-\delta)}& 0 \\
        0                   & e^{-\beta(T-\delta)}
    \end{pmatrix}=$$
    
    $$\begin{pmatrix}
        0& 1 \\
        -1   & 0
    \end{pmatrix}\begin{pmatrix}
        e^{\alpha(T-\delta)}& 0 \\
        0    & e^{-\beta(T-\delta)}
    \end{pmatrix}=\begin{pmatrix}
      0& e^{-\beta(T-\delta)}\\
          -e^{\alpha(T-\delta)}    & 0    \end{pmatrix}
$$

Now we can find the eigenvalues of the monodromy $Y(T)$ and they will be the caracteristic multipliers of (\ref{Pper})
$$Y(T)-\lambda I=\begin{pmatrix}
      -\lambda& e^{-\beta(T-\delta)}\\
          -e^{\alpha(T-\delta)}    &  -\lambda    \end{pmatrix}.$$

The caracteristic polynomial is given by $p(\lambda)\doteq \lambda ^2+e^{(\alpha-\beta)(T-\delta)}$. Since $\beta >\alpha$ this implies that 
$$|\lambda|=\sqrt{e^{(\alpha-\beta)(T-\delta)}}<1.$$
Therefore $\dot y = (A+D(t))y$ is uniformly asymptotic stable.

\section{\bf Stabilizing Unstable Linear ODE in Infinite Dimensions.}\label{unstable linear ODE}

There is a classical example in Operator Theory due to S. Kakutani of a bounded operator in an infinite-dimensional Hilbert space whose spectrum shrinks drastically from a disk to a single point under an arbitrarily small bounded perturbation. The example can be found in \cite{Rickart} (p. 282) and \cite{Halmos} (p. 248) and it is also described in \cite{8Ro-Sola}, where the present authors recently used it to build an example of the possibility of nonlinear stabilization of an unstable linear map under Fr\'{e}chet differentiability hypotheses. It is also briefly described below. The purpose of the present section is, by means of two examples, to use the ideas of Kakutani's example to show this drastic stabilization in linear ordinary differential equations in infinite dimensional Hilbert spaces, of the form
\begin{equation}\label{PK}
  \dot{x}(t)=Ax(t)+B(t)x(t),
\end{equation}
when the system $ \dot{x}(t)=Ax(t)$ is unstable and the perturbation $B(t)$ is small in some senses. Roughly speaking, we could say that the examples of this section show that while stability is a robust feature, instability does not need to be so.

Let us describe briefly the example of Kakutani with the notations and choices of \cite{8Ro-Sola}. In a real separable Hilbert space $H$ with a Hilbert orthonormal basis $(e_n)_{n\ge 1}$ a weighted shift operator $W\in\mathcal{L}(H)$ is a bounded linear operator defined by the relations $We_i=\alpha_ie_{i+1}$ for a bounded sequence of real numbers $(\alpha_n)_{n\ge 1}$. One readily sees that
\begin{equation}\label{norms}
  \|W\|=\sup\{|\alpha_n|\}\  \text{ and }  \ \|W^k\|=\sup\{|\alpha_n\alpha_{n+1}\cdots\alpha_{n+k-1}|\}.
\end{equation}
We choose first the sequence $\varepsilon_m=M/K^{m-1}$ for some $M>0$ and some $K>1$, and define a weighted shift $W_{\varepsilon}$ by $\alpha_n=\varepsilon_{m}$ if $n=2^{(m-1)}(2\ell+1)$, where $\ell$ is a non-negative integer. This sophisticated way of distributing the numbers $\varepsilon_m$ into a sequence $\alpha_n$ makes a number $\varepsilon_m$ to appear for the first time in the $\alpha_n$ sequence at the position $n=2^{(m-1)}$ and from that position onwards to appear periodically, infinitely many times, with a period of $2^m$.

Then, one also defines the weighted shifts $L_m$ by a sequence of weights $\alpha_n$ that are all of them equal to zero, except at the positions $n=2^{(m-1)}(2\ell+1)$, where $\ell$ is a non-negative integer, where $\alpha_n=\varepsilon_m$. With this choice, the operator $W_\varepsilon-L_m$ is also a weighted shift, and it has zeroes along its sequence of weights, distributed each $2^m$ places, and starting at the $2^{(m-1)}$ position. This means, according to \ref{norms}, that $W_\varepsilon-L_m$ is nilpotent of index $2^m$, $(W_\varepsilon-L_m)^{2^m}=0$. Consequently, its spectral radius $\rho(W_\varepsilon-L_m)=0$. One can also obtain, after some work, that $\rho(W_\varepsilon)=M/K$ and that the spectrum $\sigma(W_\varepsilon)$ is the whole disk of radius $M/K$ centered at zero. Concerning the norms, by using \eqref{norms} one gets that $\|W_\varepsilon\|=M$ and $\|-L_m\|=\varepsilon_m$.

In this way, Kakutani's example shows the existence of a bounded linear operator $W_\varepsilon$ with positive spectral radius that is approximated, in the operator norm, by a sequence $W_\varepsilon-L_m$ of operators whose spectrum reduces to the single point $0$.

Our fist example of translation of these ideas to \eqref{PK} is very simple. Let us choose a number $R$ and the previous numbers $M$ and $K$ in such a way that $0<R-M/K<R<1<R+M/K$ and with these choices define the new operator $T=R\,I+W_\varepsilon$, where $I$ is the identity operator. The spectrum of $T$ is a disk of radius $M/K$ centered at the point $R$. This spectrum intersects the exterior of the unit circle and lies entirely in the half-plane $\hbox{\sl Re}\, z\ge R-M/K>0$. Because of this last property, the operator $A\doteq \log(T)$ can be defined, and by the Spectral Mapping Theorem
\begin{equation}\label{Kunstable}
  \|e^{tA}\|\ge\rho(e^{tA})= e^{t\log(R+M/K)},
\end{equation}
which is unstable since $R+M/K>1$.

We construct now the sequence of operators $S_m=R\ I+W_\varepsilon-L_m$. All of these operators have their spectra reduced to the single point $z=R$, and these operators converge in the operator norm to $T=R\,I+W_\varepsilon$, which spectrum is the disk of radius $M/K$ centered at $z=R$. If we take now $A_m=\log(R\,I+W_\varepsilon-L_m)$, we again have that the sequence $A_m$ tends to $A=\log(T)$ as $m\to\infty$ in the operator norm, by the continuity of the logarithm. Also, by the properties of the exponential, perhaps by using adapted norms, for all $\delta>0$ and all $m$, there exists a number $D_{m,\delta}$ such that
\begin{equation}\label{Kstable}
  \|e^{tA_m}\|\le D_{m,\delta}\ e^{t(\log(R)+\delta)},
\end{equation}
which implies stability since $\log(R)<0$, and $\delta$ can be chosen small enough.

In this way we have perturbed an autonomous unstable system $\dot{x}(t)=Ax(t)$ to a new autonomous system $\dot{x}(t)=Ax(t)+(A_m-A)x(t)$, with a perturbation that can be taken as small as we wish in the operator norm, and the new system is asymptotically stable.

This example deserves to be commented in relation of Theorem 4 of \cite{KR} (p. 2704). According to that theorem, if an equation $\dot{x}(t)=A(t)x(t)$ exhibits an exponential dichotomy with nontrivial stable and an unstable part (which in particular means that it is unstable), then a new system $\dot{x}(t)=A(t)x(t)+B(t)x(t)$ will exhibit a similar dichotomy (which means that it is also unstable) if $\sup\{\|B(t)\|;t\in\R\}$ is sufficiently small, and if some compactness conditions are met, that are automatically satisfied in our case since $B$ does not depend on $t$. This robustness of the instability is broken in our example, since the spectrum of $A$ is a connected set that has points both in $\hbox{\sl Re}\, z<0$ and in $\hbox{\sl Re}\, z>0$, but it is not possible to divide it into two spectral sets by the vertical line $\hbox{\sl Re}\, z=0$. This is something very typical from infinite dimensional functional analysis, that cannot be expected in finite dimensions.

Our second example, also based on Kakutani's construction, starts with the same system $\dot{x}(t)=Ax(t)$ as above, with $A=\log(R\,I+W_\varepsilon)$ and $W_\varepsilon$, with the relations $0<R-M/K<R<1<R+M/K$, whose instability is expressed by the inequality \eqref{Kunstable} above. We want to add to it now a time-dependent perturbation $B(t)$, depending continuously on $t\ge 0$ such that $\sup\{\|B(t)\|;t\in [0,\infty)\}$ can be taken as small as we wish, but with the novelty that $\lim_{t\to\infty}\|B(t\|=0$. Despite of this, we want to obtain a system  $\dot{x}(t)=Ax(t)+B(t)x(t)$ that will be stable.

Let us name $B_m$ the operators $A_m-A$ considered above. Let us say again that $\|B_m\|\to 0$ as $m\to\infty$ and that the spectra $\sigma(A+B_m)=\sigma(A_m)=\{\log R\}$. Let us fix now one value of $\delta>0$ in \eqref{Kstable} such that  if we define $\omega=-\log(R)-\delta$ we still have $\omega >0$. For example, $\delta=-\frac{1}{2}\log R$. If we write $D_m$ for $D_{m,\delta}$ in \eqref{Kstable} we will have $\| e^{tA_m}\|\le D_me^{-\omega t}$. We do not expect the sequence $D_m$ to be bounded as $m\to\infty$.
Let us choose an index $m_0\ge 1$ and define 
\begin{equation}\label{B(t)}
B(t)=
\begin{cases}
B_{m_0+k},\hspace{.5cm}\hbox{ for } t_k\le t\le t_{k+1}-1,\\
(t_{k+1}-t)B_{m_0+k}+(t-t_{k+1}+1)B_{m_0+k+1} \hspace{.5cm}\hbox{ for } t_{k+1}-1\le t\le t_{k+1},
\end{cases}
\end{equation}
for an increasing sequence $t_k$ with $t_0=0$ and $t_k+1<t_{k+1}$, to be defined later. It is clear that $B(t)$ is a continuous function from $[0,\infty)$ to $\mathcal{L}(H)$. Since $\|B_m\|\to 0$ it is clear that 
$$E_{m_0}\doteq \sup\{\|B_m\|; m\ge m_0\}\to 0 \hbox{  as } m_0\to\infty.$$
Therefore, $\|B(t)\|\le E_{m_0}$ for all $t\ge 0$, and this can be made as small as we like by choosing $m_0$ sufficiently large.

In order to define the sequence $(t_k)_{k\ge 0}$ let us now bound the solutions of 
\begin{equation}\label{ABequation}
\begin{cases}
\dot{x}(t)=Ax(t)+B(t)x(t),\\
x(0)=x_0.
\end{cases}
\end{equation}
For $t$ between $t_k$ and $t_{k+1}-1$ we will have $A+B(t)=A_{m_0+k}$ and, because of \eqref{Kstable},
$$\|x(t)\|\le \|x(t_k)\|D_{m_0+k}e^{-\omega(t-t_k)}.$$

To fix ideas, let us start with $k=0$. For $0=t_0\le t\le t_1-1$ we can write
$\|x(t)\|\le D_{m_0}e^{-\omega t}\|x(0)\|$. Then, for $t_1-1\le t\le t_1$ we can broadly bound as 
$$\|x(t)\|\le e^{(t-t_1+1)(\|A\|+E_{m_0})}\|x(t_1-1)\|\le  e^{(\|A\|+E_{m_0})}\|x(t_1-1)\|,$$
and, putting the two parts together
\begin{equation}\label{first1}
  \|x(t)\|\le e^{(\|A\|+E_{m_0})}D_{m_0}e^{-\omega t}\|x(0)\|,
\end{equation}
which obviously implies the weaker bound
\begin{equation}\label{first2}
  \|x(t)\|\le e^{(\|A\|+E_{m_0})}D_{m_0}e^{-\frac{1}{2}\omega t}\|x(0)\|,
\end{equation}
both for $0\le t\le t_1$. Then, we continue with $t_1\le t\le t_2-1$, and for this range of $t$ we have $A+B(t)=A_{m_0+1}$ and
$$\|x(t)\|\le D_{m_0+1}e^{-\omega (t-t_1)}\|x(t_1)\|,$$
and, as before,
 $$\|x(t)\|\le e^{(\|A\|+E_{m_0})}D_{m_0+1}e^{-\omega (t-t_1)}\|x(t_1)\|,$$
now for the whole $t_1\le t\le t_2$. Putting this together with \eqref{first1} we get, again for $t_1\le t\le t_2$, 
$$\|x(t)\|\le e^{(\|A\|+E_{m_0})}D_{m_0+1}e^{-\omega (t-t_1)}e^{(\|A\|+E_{m_0})}D_{m_0}e^{-\omega t_1}\|x(0)\|,$$
that we can write again as 
$$\|x(t)\|\le e^{(\|A\|+E_{m_0})}D_{m_0+1}e^{-\omega (t-t_1)}e^{(\|A\|+E_{m_0})}D_{m_0}e^{-\frac{1}{2}\omega t_1}e^{-\frac{1}{2}\omega t_1}\|x(0)\|,$$
and at this point we see that we can choose $t_1$ large enough in such a way that $$e^{(\|A\|+E_{m_0})}D_{m_0+1}e^{-\frac{1}{2}\omega t_1}\le 1.$$

With this choice we get 
\begin{equation}\label{second1}
  \|x(t)\|\le e^{(\|A\|+E_{m_0})}D_{m_0}e^{-\omega (t-t_1)}e^{-\frac{1}{2}\omega t_1}\|x(0)\|,
\end{equation}
for $t_1\le t\le t_2$, which will be needed in the next interval, and also deduce, together with \eqref{first2} the weaker but more global bound
\begin{equation}\label{second2}
  \|x(t)\|\le e^{(\|A\|+E_{m_0})}D_{m_0}e^{-\frac{1}{2}\omega t}\|x(0)\|,
\end{equation}
now for all $t$ such that $0\le t\le t_2$.

Now we proceed inductively. Suppose that along the interval $t_{k-1}\le t\le t_k$, where $t_k$ is still to be chosen, we have obtained, as in \eqref{second1}, the bound 
\begin{equation}\label{kth1}
  \|x(t)\|\le e^{(\|A\|+E_{m_0})}D_{m_0}e^{-\omega (t-t_{k-1})}e^{-\frac{1}{2}\omega t_{k-1}}\|x(0)\|,
\end{equation}
for $t_{k-1}\le t\le t_k$, and the weaker inequality
\begin{equation}\label{kth2}
  \|x(t)\|\le e^{(\|A\|+E_{m_0})}D_{m_0}e^{-\frac{1}{2}\omega t}\|x(0)\|,
\end{equation}
for $0\le t\le t_k$. Then we analyze for $t_k\le t\le t_{k+1}$ and obtain that 
$$\|x(t)\|\le e^{(\|A\|+E_{m_0})}D_{m_0+k}e^{-w(t-t_k)}e^{(\|A\|+E_{m_0})}D_{m_0}e^{-\omega (t_k-t_{k-1})}e^{-\frac{1}{2}\omega t_{k-1}}\|x(0)\|.$$
Then we choose $t_k$ in such a way that 
$$e^{(\|A\|+E_{m_0})}D_{m_0+k}e^{-\frac{1}{2}\omega (t_k-t_{k-1)}}\le 1,$$
and obtain 
\begin{equation}\label{km1th1}
  \|x(t)\|\le e^{(\|A\|+E_{m_0})}D_{m_0}e^{-\omega (t-t_k)}e^{-\frac{1}{2}\omega t_k}\|x(0)\|,
\end{equation}
for $t_k\le t\le t_{k+1}$, and the weaker inequality
\begin{equation}\label{km1th2}
  \|x(t)\|\le e^{(\|A\|+E_{m_0})}D_{m_0}e^{-\frac{1}{2}\omega t}\|x(0)\|,
\end{equation}
for $0\le t\le t_{k+1}$.

With these choices of the $t_k$ one can make $k\to\infty$ and obtain the final bound 
\begin{equation}\label{kfinal}
  \|x(t)\|\le e^{(\|A\|+E_{m_0})}D_{m_0}e^{-\frac{1}{2}\omega t}\|x(0)\|,
\end{equation}
for all $t\ge 0$, that proves the exponential asymptotic stability of the solutions of \eqref{ABequation}.

\end{document}